# White Phase Intersection Control through Distributed Coordination: A Mobile Controller Paradigm in a Mixed Traffic Stream

Ramin Niroumand, Leila Hajibabai, Member, and Ali Hajbabaie, Senior Member

*Abstract*—This study presents a vehicle-level distributed coordination strategy to control a mixed traffic stream of connected automated vehicles (CAVs) and connected human-driven vehicles (CHVs) through signalized intersections. We use CAVs as mobile traffic controllers during a newly introduced "white phase", during which CAVs will negotiate the right-of-way to lead a group of CHVs while CHVs must follow their immediate front vehicle. The white phase will not be activated under low CAV penetration rates, where vehicles must wait for green signals. We have formulated this problem as a distributed mixed-integer non-linear program and developed a methodology to form an agreement among all vehicles on their trajectories and signal timing parameters. The agreement on trajectories is reached through an iterative process, where CAVs update their trajectory based on shared trajectory of other vehicles to avoid collisions and share their trajectory with other vehicles. Additionally, the agreement on signal timing parameters is formed through a voting process where the most voted feasible signal timing parameters are selected. The numerical experiments indicate that the proposed methodology can efficiently control vehicle movements at signalized intersections under various CAV market shares. The introduced white phase reduces the total delay by 3.2% to 94.06% compared to cooperative trajectory and signal optimization under different CAV market shares in our tests. In addition, our numerical results show that the proposed technique yields reductions in total delay, ranging from 40.2% - 98.9%, compared to those of a fully-actuated signal control obtained from a state-of-practice traffic signal optimization software.

## I. INTRODUCTION

Recent advancements in connected automated vehicle (CAV) technologies have the promise of significant improvements in traffic operations [1]–[3]. Obtaining online vehicle data can reduce delay at signalized intersections (e.g., [4], [5]) and signal-head-free control logic in a 100% CAV environment can maximize the intersection capacity [6], [7]. Similarly, optimally controlling the movements of CAVs through roundabouts can significantly increase their capacity [8], [9]. Since safety and policy-related challenges can hinder the transformative potential of CAVs, a smooth transition from CHVs to CAVs is critical. In fact, a rich body of literature has dealt with this transition in signalized intersections [10]–[12]. Such studies have mainly considered existing signal phases: green, yellow, and red, with various CAV market penetration rates. They either (a) optimize CAV trajectories to go through the intersection more efficiently [13], [14] or (b) jointly optimize trajectories and signal timing parameters [15], [16]. These studies have reported significant improvements in traffic operations, but they only allow one of the conflicting movements to be processed at a time regardless of the CAV market penetration rate. Besides, the existing studies mostly focus on intersections with simple layouts (e.g., one-way streets [15], [17], no turning movements [18]), or low traffic volumes [12]–[14] to tackle the complexities, particularly those of joint optimizations. Therefore, some studies form platoons of vehicles and optimize the trajectories of platoon leaders to reduce the complexity of the problem [20], [21].

There is research evidence that sufficient CAVs in traffic stream can participate in traffic control as mobile traffic controllers [22], [23] by forming platoons of vehicles and navigating them through the intersection, see Fig. 1. This paper utilizes the mobile controller paradigm and benefits from the white phase concept introduced in Niroumand et al. [22]. CAVs form platoons of CHVs, communicate with other CAVs to negotiate a safe trajectory and navigate the platoons through the intersection. CHVs only follow their immediate front vehicle to go through the intersection. This cannot happen during a green phase as CHVs may collide, which highlights the need for the newly introduced signal indication to communicate with CHVs. Note that other signal indications such as flashing green may be used to avoid the need for a change in signal heads. In-vehicle communications can be used as well. We used white phase in this paper for the ease of communication.

This study first modifies the white phase formulation presented in Niroumand et al. [22] to increase the possibility of activation of simultaneous white phases for conflicting movements that can lead to reduction in intersection total delay. Then, presents a vehicle-level distributed coordination strategy that leverages the computational power of each vehicle to address the real-time and scalability requirements of intersection traffic control. Therefore, intersections with more realistic layouts and higher traffic demand levels can be

R. Niroumand and A. Hajbabaie are with the Department of Civil, Construction, and Environmental Engineering, North Carolina State University, Raleigh, NC 27606 USA (e-mail: rniroum@ncsu.edu and ahajbab@ncsu.edu).

L. Hajibabai is with the Department of Industrial and Systems Engineering, North Carolina State University, Raleigh, NC 27606 USA (e-mail: lhajiba@ncsu.edu).
.

addressed. Each CAV receives the planned trajectory of other CAVs and the estimated trajectory of all CHVs to optimize its own trajectory. To ensure solution feasibility and moving toward system-level optimality, CAVs form agreement on their trajectories. The trajectory of each CHV is estimated using a customized car-following model. Additionally, all vehicles vote on signal timing parameters (i.e., white, green, or red) for all intersection movements based on their planned/estimated trajectories and form an agreement on signal indications to be displayed on signal heads. The vehicles' mixed-integer non-linear programs (MINLPs) are linearized and embedded into a receding horizon control to capture the dynamic nature of signal and trajectory optimizations. The proposed distributed methodology iteratively solves the problem to reach agreement on both signal indications and trajectory plans.

The remainder of the paper is organized as follows. A review of the relevant literature is presented in section II. Section III describes the details of the problem formulation. Then, the proposed solution technique is described in section IV. In section V, the case study is introduced, and numerical results are discussed. Finally, section VI provides a summary of the paper and concludes the findings.

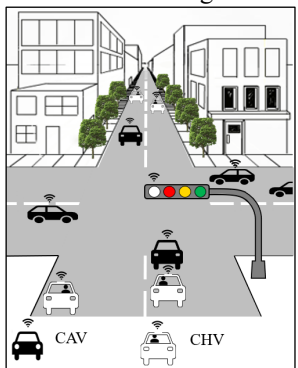

**Fig. 1.** Coordination among CAVs and CHVs.

## II. BACKGROUND

The studies on intersection control with CAVs can be categorized into three groups (a) CAV trajectory optimization, (b) signal-head-free control for 100% CAVs, and (c) joint signal timing and trajectory optimization in a mixed automated traffic stream. We have visited the first two categories in the introduction very briefly and will focus on the third group, which is the most relevant to this study.

### A. Joint trajectory and signal optimization in mixed-autonomy environments

Li and Zhou [10] have jointly optimized signal timing plans and CAV trajectories in a mixed-autonomy environment. They represented the traffic dynamics and signal timing constraints with a phase-time traffic hypernetwork model to reduce the complexity of the problem. They developed a sequential branch and bound search algorithm to solve the problem more efficiently. Their case study results indicate their proposed framework can reduce the intersection total delay by 4.2% in comparison with fixed-time signal control with a 5% CAV market share. Li et al. [24] have optimized the trajectory of CAVs and signal timing plans with the presence of conventional vehicles. They used a bi-objective and multi-stage model to minimize the traffic delay by optimizing signal timing plans and provide optimal trajectories for that signal timing plan. They utilized a hybrid heuristic method that combines genetic algorithm and particle swarm optimization. Their methodology provides optimal trajectories from each vehicle's perspective without coordination among vehicles, therefore, the solution may not improve system-level performance. Their numerical experiments on an arterial street with four intersections show 9.7%-12.2% and 0.6%-7.9% reduction in travel delay and energy consumption, respectively compared to fixed-time signal control with 10% electric vehicles.

Liang et al. [12] have developed a bi-level rolling-horizon methodology to jointly optimize signal phase and timing plan and CAV trajectories while considering traditional human-driven vehicles, connected human-driven vehicles, connected human-driven vehicles which receive speed guidance (human-cooperative vehicles), and autonomously-driven connected vehicles. The proposed methodology enumerates all possible departure sequences of naturally happening platoons to find the departure sequence with the least delay at the upper level. At the lower level, they first convert each departure sequence to a signal phase and timing and find the appropriate speed guidance for autonomous and cooperative human-driven vehicles to minimize the total number of stops. Then, the signal timings are adjusted based on the calculated speed guidance and the delay associated with the signal phase, and timing is calculated and sent back to the upper level. They tested the proposed methodology on an isolated intersection with exclusive left-turning lanes and showed that the intersection total delay and number of stops decreases as CV penetration rate increases and autonomous and cooperative human-driven vehicle penetration rate increases, respectively. Pourmehrab et al. [25] have proposed a framework to adjust an adaptive intersection control considering vehicles' arriving time to maximize the utilization of green phases. Their methodology provides a feasible signal phase and timing plan by using a set of rules to process all incoming vehicles. In addition, they split the fleet into platoons and only optimize the trajectories of automated leaders to minimize their travel delay. Their numerical experiments on a four-leg isolated intersection show a 38%-52% reduction in average travel time in comparison with a fully actuated signal control. Guo et al. [11] have proposed a two-step approach to find near-optimal signal timing plans and CAV trajectories with different CAV market penetration rates. The proposed approach finds near-optimal signal timing plans to minimize intersection delay in the first step using dynamic programming. In the second step, they generate CAV trajectories and estimate movements of CHVs based on the determined signal plans using shooting heuristic. Their case study results indicate that the proposed approach reduces average travel time and fuel consumption by up to 35.7% and 31.5%, respectively compared to adaptive signal control. A sensitivity analysis also shows that increasing the CAV penetration rate significantly contributes to improvements in intersection performance. This study can handle a mixed-autonomy traffic stream with different CAV penetration rates, however, it stops vehicles behind the red light even in a fully automated environment which leads to higher delay compared to signal-free intersection control methods. Qian et al. [23] extended a priority-based signal-free control logic to accommodate human-driven vehicles (HVs) while producing collision-free trajectories. Each CAV sends a request to the intersection controller once entered the cooperative area, in addition, infrastructure sensors are assumed to detect HVs and send a virtual request on behalf of HVs. A human-driven vehicle can enter the intersection if it is assigned the highest priority among vehicles behind the intersection stop bar or it is assigned to a virtual platoon led by a CAV. Note that admitted human driven vehicles are notified by the green signal indication. This framework can operate fully automated fleets like signal-free intersection control methods and can handle mixed-autonomy fleets by using red and green signal indications. However, using a priority-based method and forming long platoons can reduce the efficiency of the proposed framework. Moradi et al. [26] have developed an integrated platoon-based round-robin algorithm for intersections within a mixed-autonomy environment while prioritizing special vehicles. They further improved their methodology by incorporating a speed advisory mechanism for CAVs. Their methodology outperforms fixed-time and actuated intersection control methods with 5% and 35% CAV penetration rate, respectively. Du et al. [27] have proposed a joint signal and trajectory optimization framework that improves the traffic efficiency and energy saving. Their methodology optimizes the signal timing plans to minimize vehicle delay at the macro level and optimizes the CAV trajectories to minimize fuel consumption at the micro level. A comparison with the Cooperative Adaptive Cruise Control (CACC) with a fixed time signal control shows 6%-14.5% fuel consumption reduction with 20% CAV penetration rate. Tajalli et al. [28] Have developed a methodology to jointly optimize the signal timing plans and CAV trajectories in mixed autonomy environments. Their methodology decomposes the joint signal and trajectory optimization program into smaller lane-level optimization problems using a Lagrangian relaxation technique. They reported 5%-51% reduction in average travel time while having duality gaps of at most 0.1%.

Rey and Levin [29] have proposed an intersection control method for mixed-autonomy traffic streams. They have introduced a new "blue phase" during which only CAVs can access the intersection through dedicated lanes and their movements are controlled to prevent collisions with conflicting vehicles. Note that, both CAVs and HVs can enter the intersection during green phases. Their case study on a network with 25 intersections revealed that their framework can outperform existing optimized signalized intersections with at least 60% CAVs. Including CAV-dedicated lanes can dramatically decrease intersection capacity in low CAV penetration rates. Moreover, stopping all approaches to operate CAVs will require a high CAV market penetration rate to be beneficial to the intersection performance as confirmed in their case study. Niroumand et al. [22] have developed a methodology to operate a mixed-autonomy fleet through isolated signalized intersections. They have introduced a new "white phase" that is activated when there are enough CAVs in the intersection neighborhood to form vehicle-groups. During white phases, groups of CHVs led by CAVs are operated from conflicting movements simultaneously and their safety are preserved by controlling the movements of CAVs. Conventional red and green signal indications are used when there are not enough CAVs to guarantee the safety of the traffic stream by controlling the movements of CAVs. Their case studies on an isolated intersection indicate that the proposed framework works under all CAV market penetration rates. They reported a 19.6% - 96.2% reduction in intersection total delay compared to fully-actuated signal control. They further analyzed the effects of the white phase, autonomous driving behavior, and connectivity on the intersection mobility and safety performances ([30]–[32]) using the formulation and solution technique introduced in [22]. The proposed framework is flexible to operate vehicles during green phases or assign white phases to the conflicting movements and operate the intersection like a signal-free intersection. However, their case study results show that the white phase activation rate is low even in relatively high CAV penetration rates due to limited flexibility of phase transitions. In addition, they have used a central approach that does not scale well with the size of the problem. Moreover, their framework requires a central computation unit while CAVs have computational powers to be used.

### B. Summary of the literature and contribution of the paper

Some existing studies focus on signal-free intersection control logics to utilize the maximum capacity of intersections. Although signal-free control logics promise significant improvements in intersection performance measures, they require 100% CAV market share. Several studies in the literature jointly optimized the signal timing plans and CAV trajectories in mixed-autonomy traffic streams. However, they have used central methodologies to solve the joint signal and trajectory optimization program. These studies either (1) lead to reductions in intersection capacity at low CAV market penetration rates, (2) must stop vehicles behind the intersection stop bar even with 100% CAVs, and (3) do not scale well with the size of the problem.

This paper enhances the white phase formulation proposed by Niroumand et al. [22] to increase the possibility of the white phase activation and develops a vehicle-level agreement-based distributed methodology to solve it. To this end, we assumed that all vehicles are connected and have computational power to solve a simple vehicle-level optimization program. We hypothesize that the approach improves intersection capacity at low CAV market shares, reduces the need to use red phases, and finds solutions in real-time. Each CAV solves an MINLP to plan for its safe passage through the intersection and proposes a signal timing plan for the entire intersection movements. In addition, CHVs estimate their trajectories using a customized car-following model and participate in the voting process by solving an optimization problem with trajectory and signal timing variables. Unlike existing studies, this framework uses the available computational power of vehicles and solves the signal timing and trajectory optimization problem in a distributed manner. Therefore, the framework scales well with the size of the problem and can produce real-time solutions to be implemented in the real world.

## III. METHODOLOGY

This section first introduces a distributed vehicle-level MINLP for intersection control with a mixed fleet of CAVs and CHVs. Each vehicle-level model optimizes/estimates the trajectory of a single vehicle. The solutions of vehicle-level models do not necessarily lead to the system-level optimality. Therefore, subsection B introduces a coordination scheme where the system-level solutions are achieved by solving vehicle-level optimization models and coordinating their decisions via information exchange.

### A. Problem Formulation

The program aims to jointly optimize the (i) signal timing plan of each intersection and (ii) trajectory of each CAV using shared vehicle states (i.e., vehicle locations and speeds over a planning horizon) in the intersection neighborhood. Furthermore, CHVs are assumed to follow the vehicle immediately in front of them and react to signal indications.

The spatial and temporal elements of the problem are defined as follows. We let $T$ and $\hat{T}$ respectively denote the sets of discrete time steps for making trajectory and signal timing decisions (similar to [22]). Note that the size of $T$ and $\hat{T}$ are determined based on the study period $N$ and signal and trajectory updating time step lengths $\Delta T$ and $\Delta \hat{T}$. Traffic in the intersection is defined by $\mathcal{H}(L, I)$, where $L$ denotes the set of all lanes in the intersection neighborhood and $I$ represents the set of all vehicles. Additionally, $I_l \subset I$ defines the set of all vehicles on lane $l \in L$, where $I'_l$ and $\hat{I}_l$ represent CAVs and CHVs on lane $l$, respectively. We define $C_l$ as a set of all lanes in conflict with lane $l \in L$ and use $P_l$ to represent the set of all vehicle-groups on lane $l \in L$.

We define two state variables: $x_{il}^t$ represents the location of vehicle $i \in I_l$ on lane $l \in L$ at time $t \in T$ and $v_{il}^t$ denotes the speed of vehicle $i \in I_l$ on lane $l \in L$ at time $t \in T$. State variable $x_{il}^t$ identifies vehicles who pass the intersection stop bar using binary variable $\gamma_{il}^t$, where $\gamma_{il}^t = 0$ if $x_{il}^t \leq b$ or 1 otherwise (the intersection stop bar is passed) for all $i \in I_l, t \in T$. We let decision variable $a_{il}^t$ denote the acceleration rate of vehicle $i \in I_l$ on lane $l \in L$ at time $t \in T$. We estimate the movements of CHVs with a given set of signal timing parameters using a customized linear car-following model [22], [33]:

$$a_{il}^t = \max\left\{\underline{a}, \frac{\underline{v}-v_{il}^t}{\Delta T}, \min\left\{\begin{array}{c}\overline{a}, \\ \frac{\overline{v}-v_{il}^t}{\Delta T}, \\ \alpha_1(v_{il}^t - \hat{v}_{i-1,l}^t) + \alpha_2((x_{il}^t - \hat{x}_{i-1,l}^t - \mathcal{L}) - D - v_{il}^t \hat{\tau}), \\ \alpha_1\left((w_l^n + g_l^n)\overline{v} - v_{il}^t\right) + \alpha_2\left((b_l - x_{il}^t) - (1 - w_l^n - g_l^n)S\right) + M\gamma_{il}^t,\end{array}\right\}\right\},$$

$$\forall i \in \hat{I}_l, l \in L, t \in T, n \in \hat{T}, \quad (1)$$

where $M\gamma_{il}^t$ eliminates the connection between signal and vehicles when they pass the intersection stop bar.

We also define decision variables $g_l^n = \{0,1\}$ and $w_l^n = \{0,1\}$ to represent the green and white signal timing decisions for all $l \in L$ and $n \in \hat{T}$, respectively. TABLE I defines the remaining variables and parameters used in the proposed formulation. Fig. 2 shows some defined notations in the layout of an isolated intersection.

TABLE I NOTATIONS USED IN THE PROPOSED FORMULATION.

| | |
|---|---|
| **Sets** | |
| $L$ | Set of all lanes |
| $T$ | Set of all time steps for trajectory decisions ($T = \{1,2,...,N/\Delta T\}$) |
| $\hat{T}$ | Set of all time steps for signal timing decisions ($\hat{T} = \{1,2,...,N/\Delta\hat{T}\}$) |
| $I_l$ | Set of all vehicles on lane $l \in L$ |
| $I'_l$ | Set of all CAVs on lane $l \in L$ |
| $\hat{I}_l$ | Set of all CHVs on lane $l \in L$ |
| $C_l$ | Set of all conflicting lanes with lane $l \in L$ |
| $P_l$ | Set of vehicle-groups on lane $l \in L$ |
| **Decision Variables** | |
| $a_{il}^t$ | Acceleration of vehicle $i \in I_l$ on lane $l \in L$ at time $t \in T$ |
| $g_l^n$ | Binary green signal status; 1 if signal for lane $l \in L$ at time $n \in \hat{T}$ is green or 0 otherwise |
| $w_l^n$ | Binary white signal status; 1 if signal for lane $l \in L$ at time $n \in \hat{T}$ is white or 0 otherwise |
| $y_l^n$ | Binary white signal status; 1 if signal for lane $l \in L$ at time $n \in \hat{T}$ is yellow or 0 otherwise |
| **Variables** | |
| $x_{il}^t$ | The location of vehicle $i \in I_l$ on lane $l \in L$ at time $t \in T$ |
| $v_{il}^t$ | The speed of vehicle $i \in I_l$ on lane $l \in L$ at time $t \in T$ |
| $\gamma_{il}^t$ | 1 if vehicle $i \in I_l$ on lane $l \in L$ has passed intersection by $t \in T$ or 0 otherwise |
| $h_{ql}^t$ | The head location of vehicle-group $q \in P_l$ on lane $l \in L$ at time $t \in T$ |
| $e_{ql}^t$ | The tail location of vehicle-group $q \in P_l$ on lane $l \in L$ at time $t \in T$ |
| $\zeta_{ql}^t$ | The length of vehicle-group $q \in P_l$ on lane $l \in L$ at time $t \in T$ |
| **Parameters** | |
| $\underline{a}, \overline{a}$ | Lower and upper bounds for acceleration |
| $\underline{v}, \overline{v}$ | Minimum and maximum speeds |
| $\mathcal{L}$ | Fixed vehicle length |
| $D$ | Minimum safety distance between two consecutive vehicles on the same lane |
| $S$ | Minimum distance between vehicles and the intersection stop bar |
| $\rho$ | Minimum safety distance between vehicle-groups of conflicting lanes |
| $\hat{\tau}$ | Reaction time for CHVs |
| $\tau'$ | Reaction time for CAVs |
| $F_{kk'}$ | The location of conflict point between lane $k \in L$ and lane $k' \in C_k$ based on the coordination of lane $k$ |
| $b$ | The location of intersection (i.e., the distance of intersection *stop* bar from origin) |
| $\overline{G}$ | Maximum green time |
| $\underline{G}$ | Minimum active time (green plus subsequent white phase) |
| $\underline{W}$ | Minimum white time |
| $Y$ | Yellow time duration |
| $\Delta T$ | Time step size for trajectory decisions |
| $\Delta\hat{T}$ | Time step size for signal decisions |
| $r_l$ | The destination location of vehicles on lane $l \in L$ |
| $\alpha_1, \alpha_2$ | Car-following parameters |
| $\mathcal{R}$ | Duration of all-red time |
| $\xi_{ql}$ | The order of the first vehicle of vehicle-group $q \in P_l$ on lane $l \in L$ |
| $\mu_{ql}$ | The order of the last vehicle of vehicle-group $q \in P_l$ on lane $l \in L$ |
| $\beta_{il}^t$ | 1 if vehicle $i \in I_l$ on lane $l \in L$ is a member of a vehicle-group or 0 otherwise |
| $\hat{\gamma}_{jl}^t$ | Shared input parameter: 1 if vehicle $j \in I_l$ on lane $l \in L$ has passed intersection by $t \in T$ or 0 otherwise |
| $\hat{x}_{jl}^t$ | Location of vehicle $j \in I_l$ on lane $l \in L$ as shared input parameter |
| $\hat{g}_{ljk}^n$ | Green signal indications on lane $l \in L$ at time $n \in \hat{T}$ voted by vehicle $j \in I_l$ on lane $k \in L$ |
| $\hat{w}_{ljk}^n$ | White signal indications on lane $l \in L$ at time $n \in \hat{T}$ voted by vehicle $j \in I_l$ on lane $k \in L$ |
| $\mathcal{J}$ | Iteration counter for agreement process |
| $\overline{\zeta}$ | The maximum length of vehicle-groups |

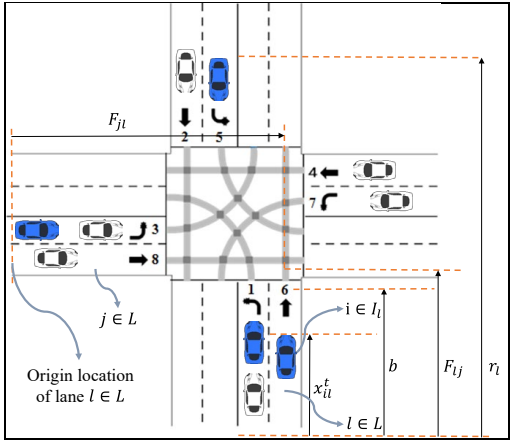

**Fig. 2.** An isolated intersection layout including some of the notations

The problem formulation follows. The objective is to minimize the distance between the location of human-driven vehicle $i \in \hat{I}_l$ (or automated vehicle $i \in I'_l$) and its pre-defined destination. A term $\omega|v_{il}^{t+1} - v_{il}^t|$ is added to the objective function $Z'_{il}$ of CAVs to minimize the speed variation of vehicle $i \in I'_l$ for driving comfort purposes. Parameter $\omega$ is the speed-location conversion factor.

$$\hat{Z}_{il} = \min_{a,g,w} \sum_{t \in T}(r_l - x_{il}^t), \qquad \forall i \in \hat{I}_l, l \in L, \quad (2)$$

$$Z'_{il} = \min_{a,g,w} \sum_{t \in T}(r_l - x_{il}^t) + \omega|v_{il}^{t+1} - v_{il}^t|, \quad (3)$$
$$\forall i \in I'_l, l \in L.$$

Subject to (1) and:

$$x_{il}^{t+\Delta T} = x_{il}^t + v_{il}^t \Delta T + \frac{1}{2}a_{il}^t \Delta T^2, \quad (4)$$
$$\forall i \in I_l, l \in L, t \in T,$$

$$v_{il}^{t+\Delta T} = v_{il}^t + a_{il}^t \Delta T, \qquad \forall i \in I_l, l \in L, t \in T, \quad (5)$$

$$\underline{a} \leq a_{il}^t \leq \overline{a}, \qquad \forall i \in I'_l, l \in L, t \in T, \quad (6)$$

$$\underline{v} \leq v_{il}^t \leq \overline{v}, \qquad \forall i \in I'_l, l \in L, t \in T, \quad (7)$$

$$\hat{x}_{i-1,l}^t - x_{il}^t \geq D + \mathcal{L} + \tau' v_{il}^t, \qquad \forall i \in I'_l, l \in L, t \in T, \quad (8)$$

$$b - x_{il}^t \geq S + \Delta T v_{il}^t - (g_l^n + w_l^n)M - \gamma_{il}^t M, \quad (9)$$
$$\forall i \in I'_l, l \in L, t \in T, n \in \hat{T},$$

$$h_{qj}^t = \hat{x}_{\xi_{qj}j}^t, \quad (10)$$
$$\forall q \in P_j, j \in C_l, l \in L, t \in T,$$

$$e_{qj}^t = \hat{x}_{\mu_{qj}j}^t - \mathcal{L}, \quad (11)$$
$$\forall q \in P_j, j \in C_l, l \in L, t \in T,$$

$$\zeta_{qj}^t = h_{qj}^t - e_{qj}^t, \quad (12)$$
$$\forall q \in P_j, j \in C_l, l \in L, t \in T,$$

$$|h_{qj}^t - F_{jl}| + |e_{qj}^t - F_{jl}| + |x_{il}^t - F_{lj}| + |x_{il}^t - \mathcal{L} - F_{lj}| \geq \zeta_{qj}^t + \mathcal{L} + 2\rho + M(\gamma_{il}^t + \hat{\gamma}_{\xi_{qj}j}^t + w_l^n + w_j^n - 4), \quad (13)$$
$$\forall i \in I'_l, l \in L, q \in P_j, j \in C_l, t \in T, n \in \hat{T},$$

$$g_l^n + g_{l'}^n + w_{l'}^n \leq 1, \quad (14)$$
$$\forall l \in L, l' \in C_l, n \in \hat{T},$$

$$\sum_{z=n+1}^{n+\underline{G}} g_l^z + w_l^z \geq (g_l^{n+1} - g_l^n)\underline{G}, \quad (15)$$
$$\forall l \in L, n \in \hat{T},$$

$$\sum_{z=n+1}^{n+\underline{W}} w_l^z \geq (w_l^{n+1} - g_l^n - w_l^n)\underline{W}, \quad (16)$$
$$\forall l \in L, n \in \hat{T},$$

$$\sum_{z=n}^{n+\overline{G}} g_l^z \leq \overline{G}, \qquad \forall l \in L, n \in \hat{T}, \quad (17)$$

$$\sum_{z=n}^{n+Y} y_l^z \leq Y, \qquad \forall l \in L, n \in \hat{T}, \quad (18)$$

$$\sum_{z=n+1}^{n+Y} y_l^z \geq (g_l^n - g_l^{n+1} - w_l^{n+1})Y, \quad (19)$$
$$\forall l \in L, n \in \hat{T},$$

$$\sum_{z=n+1}^{n+Y} y_l^z \geq (w_l^n - w_l^{n+1} - g_l^{n+1})Y, \quad (20)$$
$$\forall l \in L, n \in \hat{T},$$

$$\sum_{z=n}^{n+\mathcal{R}-1}(g_l^z + g_{l'}^z + w_l^z + w_{l'}^z) - 2\mathcal{R}y_l^n + 2\mathcal{R}y_l^{n-1} \leq 2\mathcal{R}, \quad (21)$$
$$\forall l \in L, l' \in C_l, n \in \hat{T},$$

$$\sum_{z=n}^{n+\mathcal{R}-1}(w_{l'}^z) - 2\mathcal{R}g_l^n + 2\mathcal{R}g_l^{n-1} \leq 2\mathcal{R}, \quad (22)$$
$$\forall l \in L, l' \in C_l, n \in \hat{T},$$

$$w_l^n \leq 1 - (\hat{\gamma}_{j-1,l}^t - \hat{\gamma}_{jl}^t) + w_l^{n-1}, \quad (23)$$
$$\forall j \in \hat{I}_l, l \in L, n \in \hat{T}, t \in T,$$

$$w_l^n \leq 1 - (\hat{\gamma}_{j-1,l}^t - \hat{\gamma}_{jl}^t) + \beta_{jl}^t, \quad (24)$$
$$\forall j \in \hat{I}_l, l \in L, n \in \hat{T}, t \in T,$$

Constraints (4) and (5) define fundamental motion equations to update the location $x_{il}^{t+\Delta T}$ and speed $v_{il}^{t+\Delta T}$ of vehicle $i \in I_l$ at time $t + \Delta T$ based on acceleration rate $a_{il}^t$ at time $t$. Constraints (6) and (7) ensure that CAVs are enforced to take their acceleration and speed within the threshold of pre-defined minimum and maximum values. Constraints (8) ensure the safety of vehicle $i \in I'_l$ on lane $l \in L$ at time $t \in T$ based on its distance from preceding vehicle $i - 1 \in I_l$, where $D$ is a constant safety distance between two consecutive vehicles on the same lane, $L_v$ is the vehicle length, and $\tau' v_{il}^t$ is the distance traveled by the CAV during reaction time $\tau'$, given the location of vehicle $i - 1 \in I_l$ on lane $l$ at time $t$ ($\hat{x}_{i-1,l}^t$) that is fed to the optimization program as an input parameter.

Constraints (9) will force CAVs to stop behind the intersection stop bar, keeping a distance $S$ when the signal indication is red. Note that vehicles are allowed to enter the intersection conflicting area during green and white phases. Since the formulation uses discrete time steps, a vehicle can jump from one side of the intersection stop bar to the other side during one time step when the signal is red, therefore, we add the term $\Delta T v_{il}^t$ to prevent this jump. Besides, the term $-\gamma_{il}^t M$ ensures that CAVs who have already passed the intersection stop bar will not react to any signal indication. Note that $S$ can be as small as zero. Constraints (10) – (12) define a group of CHVs led by a CAV as a vehicle-group that starts with the CAV and ends with the last CHV in the group, where $h_{qj}^t$, $e_{qj}^t$, and $\zeta_{qj}^t$ respectively denote the head location, tail location, and length of vehicle-group $q \in P_j$ on lane $j \in C_l$ at time $t \in T$. Each vehicle group starts with a CAV as the group leader and continues until it reaches the maximum predefined vehicle group length $\bar{\zeta}$ or the next CAV in that lane. Note that each vehicle-group contains only one CAV as its leader. For instance, vehicle 1 shown in Fig. 3 is a vehicle group by itself since the next vehicle is also a CAV. The next vehicle group starts with vehicle 2 as the vehicle group leader and has vehicle 3 as its follower. The third vehicle group starts with vehicle 4 and has vehicle 5 as its follower, however, vehicle 6 is not a member of the third vehicle group since the third vehicle group has reached its maximum predefined length before reaching vehicle 6. Finally, vehicle 7 is a vehicle group by itself since it is a CAV and there is no vehicle after vehicle 7. These constraints aim to connect the vehicle-level variables (i.e., $x_{il}^t$) and vehicle-group-level variables (i.e., $h_{qj}^t$, $e_{qj}^t$, and $\zeta_{qj}^t$). Constraints (13) enforce collision-free movements based on the topology of conflicting vehicle-groups in the intersection zone when the signal is white [22]. The additional term $M(\gamma_{il}^t + \hat{\gamma}_{\xi_{qjj}}^t)$ relaxes the constraints when at least one of the vehicle-groups has not entered the conflicting area.

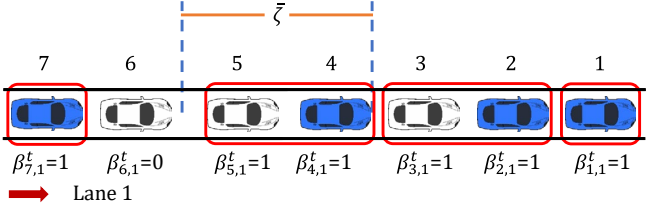

**Fig. 3.** Vehicle group formation

Additionally, constraints (14) ensure that green-white and green-green signals will not activate simultaneously for conflicting lane-groups (see Fig. 4). Note that white-white is safe for conflicting lane-groups. However, green and white signal indications can be used for non-conflicting movements simultaneously since they don't share any conflicting points.

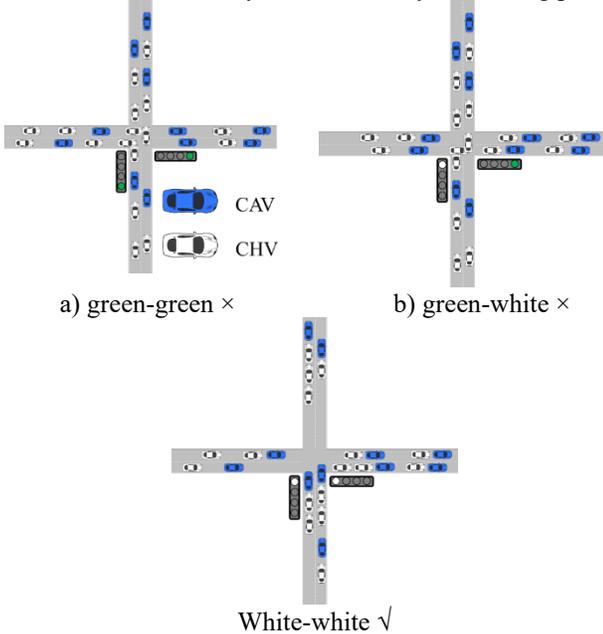

**Fig. 4.** Different signal combinations for conflicting lane groups

We introduced a set of constraints to impose a lower bound on the minimum green time in the central formulation. This set of constraints state that if a lane group starts receiving green signal indication, it should receive green signal indication for a minimum period $\underline{G}$. However, these constraints can hinder the activation of white phases, which is not beneficial in terms of intersection delay. The first vehicle behind the intersection stop bar of east bound movement shown in Fig. 5.a is a CHV, therefore, it is not possible to give white indication to east bound movement. The only option to clear the east bound queue is to give green indication to that movement and red indication to the north bound movement. Based on the minimum green constraints, this signal plan should be active for at least $\underline{G}$ minutes. However, if we switch the signal indication from green to white after operating the first vehicle behind the east bound stop bar, we can operate both east bound and north bound movements under simultaneous white phases as shown in Fig. 5.b. Hence, we replaced the minimum green time constraints by minimum active time constraints presented as constraints (15) which imposes a lower bound on the active time (green and following white times) of a lane-group when its signal indication switches from red to green. Similarly, constraints (16) impose a minimum duration for the white phase for a lane-group when its signal switches from red to white. constraints (17) define an upper bound for the green signal duration of a lane group.

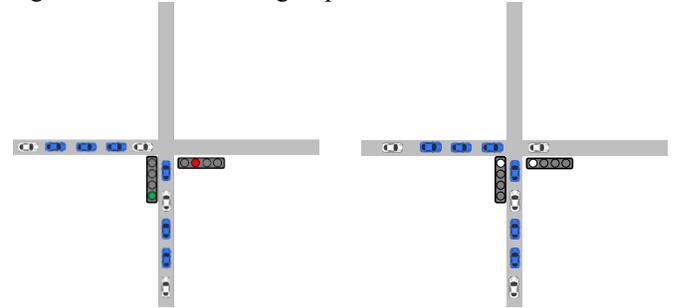

a) Minimum green time    b) Minimum active time
**Fig. 5.** Minimum green and active times

Constraints (18) define the yellow time duration while constraints (19) and (20) ensure that signal turns from green to yellow at the end of green and from white to yellow at the end of the white time, respectively. Constraints (21) impose all-red indication $\mathcal{R}$ when a lane-group switches from yellow to red. Constraints (22) are introduced to clear the intersection conflicting area while switching the signal from green to white. In particular, if the signal is green for lane $l \in L$ at time $n \in \hat{T}$ and the lane is to receive white signal at time $n + 1 \in \hat{T}$, constraints (22) will prevent conflicting lanes with lane $l \in L$ from receiving white signal for $\mathcal{R}$ duration. Therefore, the last vehicle that entered the intersection conflicting area will clear the intersection before letting the conflicting vehicle-groups enter the conflicting area during white phases.

Constraints (23) – (24) define the initiation and termination of the white phase. The phase can initiate only when the first vehicle behind the intersection stop bar is a CAV. Additionally, the white phase will be terminated for lane $l \in L$ if the first vehicle behind the intersection stop bar is not a member of a vehicle-group. For instance, the third vehicle shown in Fig. 6 is the first vehicle behind the intersection stop bar and is not a member of the vehicle-group, therefore, its lane-group cannot receive white.

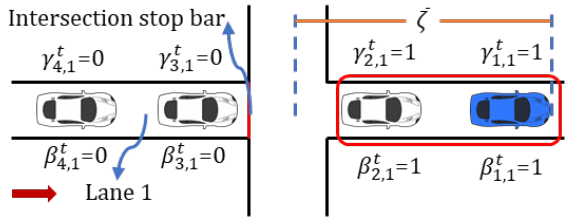

**Fig. 6.** The relative location of vehicles in reference to intersection stop bar and vehicle-groups

In summary, CAVs optimize objective function (3) subject to constraints (4) - (24). On the other hand, CHVs estimate their movements using the car-following model (1) with a given signal, or optimize objective function (2) subject to constraints(1), (4), (5), and (14) - (24) to estimate their trajectories while finding the optimal signal timing plans for themselves.

### B. Solution Technique

The proposed vehicle-level MINLP in (1) contains non-linear terms in the objective function (3) and constraints (1and (13); and seeks agreement on trajectories among vehicles with conflicting movements. This subsection first linearizes the aforementioned non-linear terms and then introduces a receding horizon structure over a planning horizon $\hat{N}$ to capture the dynamics of distributed trajectory and signal timing plans. We first define auxiliary non-negative variables $\lambda^{+t}_{il}$ and $\lambda^{-t}_{il}$ to linearize objective function (3) as follows.

$$\tilde{Z}'_{il} = \min_{a,g,w} \sum_{t \in T}(r_l - x^t_{il}) + \omega(\lambda^{+t+1}_{il} + \lambda^{-t+1}_{il}), \quad (25)$$
$$\forall i \in I'_l, l \in L,$$

$$v^{t+1}_{il} - v^t_{il} = \lambda^{+t+1}_{il} - \lambda^{-t+1}_{il}, \quad (26)$$
$$\forall i \in I'_l, l \in L, t \in T,$$

$$\lambda^{+t+1}_{il}, \lambda^{-t+1}_{il} \geq 0, \quad (27)$$
$$\forall i \in I'_l, l \in L, t \in T,$$

where (25) ensures that one of the auxiliary variables gets the value of zero and the other equals the absolute value function. Similar to [22], car-following constraints (1) will also be linearized. Additionally, we replace the safety constraints (13) by:

$$\varphi^{+tq}_{jl} + \varphi^{-tq}_{jl} + \phi^{+tq}_{jl} + \phi^{-tq}_{jl} + \varphi^{+ti}_{lj} + \varphi^{-ti}_{lj} +$$
$$\phi^{+ti}_{lj} + \phi^{-ti}_{lj} \geq \zeta^t_{qj} + \mathcal{L} + 2\rho + M\left(\gamma^t_{il} + \hat{\gamma}^t_{\xi_{qj}j} + w^n_l + w^n_j - 4\right), \quad (28)$$
$$\forall i \in I'_l, l \in L, q \in P_j, j \in C_j, t \in T, n \in \hat{T},$$

where $\varphi$ and $\phi$ are auxiliary non-negative variables that help linearize the absolute value functions corresponding to the head location and tail location of vehicle-groups, respectively.

Since the proposed optimization program contains binary variables, it is possible to provide some cuts to expedite the solution process of the optimization program. For instance, binary variable $\gamma^t_{il}$ takes the value of one if vehicle $i \in I_l$ on lane $l \in L$ has passed the intersection stop bar by time $t \in T$ and takes the value of zero otherwise. Based on the definition of this variable, we can provide following cuts to reduce the computation time of the problem.

$$\gamma^{t'}_{il} \leq \gamma^t_{il}, \quad \forall i \in I_l, l \in L, t \in T, t' < t, \quad (29)$$

$$\gamma^t_{i'l} \leq \gamma^t_{il}, \quad \forall i \in I_l, l \in L, t \in T, i' > i, \quad (30)$$

$$\gamma^t_{il} = 1, \quad if\ x^0_{il} \geq b \quad \forall i \in I_l, l \in L, t \in T, \quad (31)$$

Constraints (29) state that if vehicle $i \in I_l$ on lane $l \in L$ has not passed the intersection location by time $t \in T$ (i.e., $\gamma^t_{il} = 0$), the state variable $\gamma$ for that vehicle should take the value of zero for time period $t' \in [0, t]$. Similarly, constraints (30) indicates that if vehicle $i \in I_l$ on lane $l \in L$ has not passed the intersection location by time $t \in T$, the state variable $\gamma$ for its following vehicles should take the value of zero at time $t \in T$. Finally, constraints (30) states that if a vehicle has already passed the intersection location (i.e., $x^0_{il} \geq b$) the state variable $\gamma$ for that vehicle should take the value of one for the entire planning horizon. Note that $x^0_{il}$ is the current location of the vehicle $i \in I_l$ on lane $l \in L$ which is a parameter. Although these constraints seem primitive, they can reduce the computation time by pruning unnecessary branches in the branch and cut process used by CPLEX. Similar cuts are also introduced for the binary variables used for linearizing the safety constraints (13).

*1) Receding horizon technique*

Similar to our previous studies ([34]–[36]), a receding horizon control is implemented to solve the proposed formulation over a finite planning horizon $\hat{N}$. We first initialize the state variables $x^{t_0}_{il}$ and $v^{t_0}_{il}$ of vehicle $i \in I_l$ on lane $l \in L$ at time $t_0 \in T$. We also capture the current signal status, i.e., green $g^{n_0}_l$ and white $w^{n_0}_l$ on lane $l \in L$ at time $n_0 \in \hat{T}$. The algorithm proceeds iteratively to form agreement among all vehicles on trajectories as well as signal timing plans. Since the update time steps for trajectory and signal timing plans are different (following sets $T$ and $\hat{T}$, respectively), the algorithm reaches agreement on either (i) trajectories with given signal timing plans at time $t \in [t_0, t_0 + \hat{N}/\Delta T]$ or (ii) both signal timing plans and trajectories at times $n \in [n_0, n_0 + \hat{N}/\Delta \hat{T}]$ and $t \in [t_0, t_0 + \hat{N}/\Delta T]$. Post agreement, the selected actions in the first time step that are (1) acceleration rate $a^t_{il}$ of CAV $i \in I'_l$ on lane $l \in L$ at time $t_0 + 1$ and (2) signal indications $w^n_l$ and $g^n_l$ on lane $l \in L$ at time $n_0 + 1$ will be implemented. Note that the trajectories of CHVs are updated by Vissim and the car-following model (1) is only used to estimate the trajectories of the CHVs to be used in the trajectory optimization of CAVs. Tis mismatch between the estimated and implemented trajectories for CHVs results in a stochastic error which is captured using the receding horizon framework. After updating the acceleration rates and signal indications, the same process is repeated to reach agreement on signal timing plans and CAV trajectories at times $n \in [n_0 + 1, n_0 + 1 + \hat{N}/\Delta \hat{T}]$ and $t \in [t_0 + 1, t_0 + 1 + \hat{N}/\Delta T]$. The same process is repeated until the end of the study period. The proposed framework is detailed in Fig. 7.

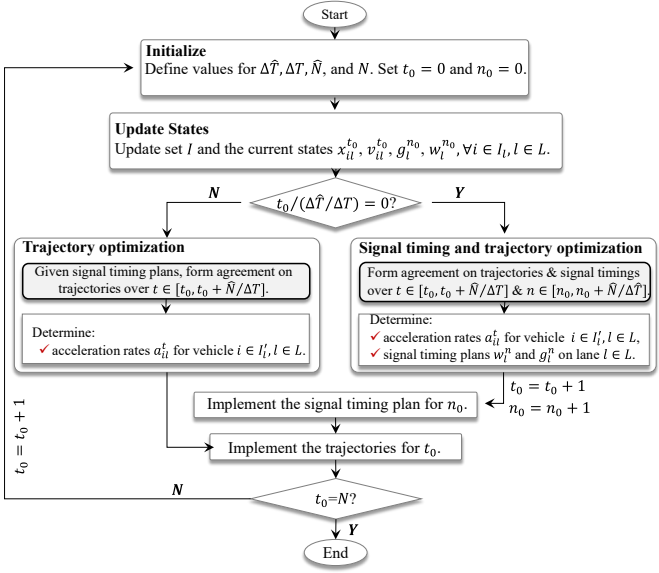

**Fig. 7.** Receding horizon framework.

*2) Agreement*

Each vehicle optimizes/estimates its own trajectory and proposes a signal timing plan for all intersection lane-groups based on the shared information obtained from the rest of the vehicles in the intersection neighborhood. The decision variables of the vehicle-level optimization models are the acceleration rate of the vehicle $a$ and green $g$ and white $w$ signal statuses for the lane groups. As a result, each vehicle chooses the signal status to maximize its own objective function. The value of the signal timing variables optimized by each vehicle is considered as its vote on signals. If the signals proposed by vehicles do not change in two consecutive iterations, an optimization model is solved to find the most voted feasible signal timing parameters to be implemented. Note that there is no need for CHVs to be assigned to a vehicle-group to be able to vote on signal timing parameters. Vehicle movements are subject to safety constraints (28) during white phases to avoid conflicting vehicle-groups. Therefore, vehicles tend to vote for a green signal for their lane-groups. Hence, vehicle-level objective functions incentivize voting for white signals. We also relax constraints (28) by introducing slack variable $\delta_{il}^{tqj}$ in the CAVs' objective function with a big coefficient $\psi$ to promote high coordination levels among CAVs. Therefore, objective functions $\hat{Z}_{il}$ and $\tilde{Z}'_{il}$ are rewritten as follows.

$$\widehat{Z^c}_{il} = \min_{a,g,w} \sum_{t \in T}(r_l - x_{il}^t) + M(a_{il}^t - d_{il}^t) - \mu \sum_{n \in \hat{T}} \sum_{j \in L} w_j^n, \quad \forall i \in \hat{I}_l, l \in L, \quad (32)$$

$$\widetilde{Z^c}'_{il} = \min_{a,g,w} \sum_{t \in T}(r_l - x_{il}^t) + \omega(\lambda_{il}^{+t+1} + \lambda_{il}^{-t+1})$$
$$-\mu \sum_{n \in \hat{T}} \sum_{j \in L} w_j^n + \psi \sum_{t \in T} \sum_{q \in P_j} \sum_{j \in C_l} \delta_{il}^{tqj}, \quad (33)$$
$$\forall i \in I'_l, l \in L,$$

Similarly, safety constraints (28) are reformulated as:
$$\varphi^{+tq}_{jl} + \varphi^{-tq}_{jl} + \phi^{+tq}_{jl} + \phi^{-tq}_{jl} + \varphi^{+ti}_{lj} + \varphi^{-ti}_{lj} + \phi^{+ti}_{lj} + \phi^{-ti}_{lj} + \delta_{il}^{tqj} \geq \zeta_{qj}^t + \mathcal{L} + 2\rho + M\left(\gamma_{il}^t + \hat{\gamma}_{\xi qj}^t + w_l^n + w_j^n - 4\right), \quad (34)$$
$$\forall i \in I'_l, l \in L, q \in P_j, j \in C_j, t \in T, n \in \hat{T},$$

We implement the following steps to reach agreement among vehicles on trajectories and signal timing plans.

**Step i** Update CAV trajectories

Given the signal timing plans, CAVs solve problem $z$ with objective (33) and constraints (4)-(12), (14)-(24), (26)-(27), and (34) to optimize their trajectories at each iteration $\mathcal{T}$. The trajectory of vehicle $i \in I'_l$ on lane $l \in L$ at time $t \in T$ at iteration $\mathcal{T} + 1$ will be updated by averaging ([37]) the trajectory from iteration $\mathcal{T}$ (i.e., $\hat{x}_{il}^{t\mathcal{T}}$) and the optimized trajectory $x_{il}^t$ to push the trajectories toward system-level agreement after a finite number of iterations $\mathcal{T}$. The maximum allowable value for the slack variable $\boldsymbol{\delta}$ is reduced at each iteration to ensure feasible trajectories.

$$\hat{x}_{il}^{t,\mathcal{T}+1} = \left(1 - \frac{1}{\mathcal{T}}\right)\hat{x}_{il}^{t\mathcal{T}} + \left(\frac{1}{\mathcal{T}}\right)x_{il}^t | z, \quad (35)$$
$$\forall i \in I'_l, l \in L, t \in T.$$

**Step ii** Update CHV trajectories

CHVs use car-following model (1) to estimate their trajectories based on given signal timing plans and shared information from the rest of the fleet. Note that, the trajectories of CHVs are not averaged since human drivers do not use the exchanged information among vehicles and only react to the implemented trajectories of CAVs. Therefore, the trajectory of vehicle $i \in \hat{I}_l$ on lane $l \in L$ at time $t \in T$ at iteration $\mathcal{T} + 1$ is updated only based on the obtained trajectories as follows.

$$\hat{x}_{il}^{t,\mathcal{T}+1} = x_{il}^t | z', \quad \forall i \in \hat{I}_l, l \in L, t \in T. \quad (36)$$

**Step iii** Share trajectory data

All vehicles share their updated trajectory $\hat{x}_{jl}^{t,\mathcal{T}+1}$ with the rest of the vehicles in the intersection neighborhood to be used at iteration $\mathcal{T} + 1$.

**Step iv** Subroutine for trajectory agreement

A trajectory agreement is reached when the trajectory changes at two consecutive iterations do not exceed a predefined value $\varepsilon$; i.e.,
$$|\hat{x}_{il}^{t,\mathcal{T}+1} - \hat{x}_{il}^{t,\mathcal{T}}| \leq \varepsilon, \quad \forall i \in I_l, l \in L, t \in T, \quad (37)$$

**Step v** Subroutine for signal agreement

At time steps $t$ where joint optimization is involved; each vehicle optimizes its own trajectory and proposes a signal timing plan for all lane-groups. To this end, CAVs solve problem $z$ while CHVs solve problem $z'$ with objective (32) and constraints (1), (4)-(5), and (14)-(24). Note that vehicles that have arrived the intersection earlier must clear the conflicting area before their following vehicles can access the intersection. Therefore, we introduce constraints (38) and (39) to make sure that vehicles respect to the votes proposed by their preceding vehicles. Constraints (38) and (39) ensure that vehicle $i \in I_l$ on lane $l \in L$ accepts the signal timing plan proposed by its preceding vehicles on the same lane; this guarantees that preceding vehicles have cleared the intersection by the time vehicle $i \in I_l$ enters the intersection. Note that constraints (38) and (39) are defining priorities only for the vehicles on the same lane group while the negotiation for signal timing parameters happens between conflicting lane groups. Green and white signal indications on lane $j \in L$ at time $n \in \hat{T}$ voted by vehicle $i \in I_l$ on lane $l \in L$ is denoted by $\hat{g}_{ljl}^n$ and $\hat{w}_{ljl}^n$. The term $1 - \hat{\gamma}_{jl}^t$ will relax the two constraints after vehicle $j$ passes the intersection stop bar. For instance, assume that vehicle 3 shown in Fig. 6 proposes green signal indication for lane 1 from time $n = 0$ to $n = 10$. In addition, assume that vehicle 3 is planning to cross the intersection stop bar by the time $t = 8$ which leads to $\hat{\gamma}_{3,1}^8 = 1$. Please note that signal time step $n = 2$ starts at the same time as the trajectory time step $t = 8$ since we update trajectories and signal indications every 0.5 and 2 seconds, respectively. Constraints (38) and (39) ensure that vehicle 4 proposes the same signal indications that vehicle 3 proposed for lane 1 to ensure that vehicle 3 has cleared the intersection by the time that vehicle 4 enters the intersection. Furthermore, the term $(1 - \hat{\gamma}_{jl}^t)$ is multiplied to the right-hand side of these constraints to relax them after the preceding vehicle has passed the intersection stop bar. Therefore, vehicle 2 should propose green signal indication for lane 1 from time step $n = 0$ to time step $n = 2$.

$$w_l^n \geq \widehat{w}_{ljl}^n (1 - \hat{\gamma}_{jl}^t), \quad j < i, l \in L, n \in \hat{T}, t \in T, \quad (38)$$

$$g_l^n \geq \hat{g}_{ljl}^n (1 - \hat{\gamma}_{jl}^t), \quad \forall j < i, l \in L, n \in \hat{T}, t \in T. \quad (39)$$

If the voted signal timing plans remain unchanged in consecutive iterations, each vehicle will solve the following problem to optimize the signal timings. The signal plans are incorporated into the optimization program of vehicles as input parameters for the next iteration. The objective is to find the signal status for each lane-group with the highest vote from vehicles on the same lane-group.

$$\mathbb{Z} = \min_{g, w} \sum_{n \in \hat{T}} \sum_{l \in L} \sum_{i \in I_l} (\mathcal{D}_{il} + 1)(1 - \hat{\gamma}_{il}^t)(|g_l^n - \hat{g}_{lil}^{n\mathcal{T}}| + |w_l^n - \widehat{w}_{lil}^{n\mathcal{T}}|) \quad (40)$$
$$s.t. (14)-(24), (38)-(39)$$

where $\mathcal{D}_{il}$ is the delay of vehicle $i \in I_l$ on lane $l \in L$, obtained from Vissim [38] as an input parameter, which imposes greater weight to vehicles with higher delay. For instance, if the delay of vehicle $i \in I_l$ on lane $l \in L$ is 1, its vote will be multiplied by $\mathcal{D}_{il} + 1 = 2$. Additionally, $\hat{g}_{lil}^{n\mathcal{T}}$ and $\widehat{w}_{lil}^{n\mathcal{T}}$ respectively denote green and white signal indications on lane $l \in L$ at time $n \in \hat{T}$ voted by vehicle $i \in I_l$ on lane $l \in L$ at iteration $\mathcal{T}$. Multiplying $1 - \hat{\gamma}_{il}^t$ captures vehicle $i$'s vote as far as it has not passed the intersection stop bar.

We utilize a linearized version of $\mathbb{Z}$ and associated constraints as follows.

$$\widetilde{\mathbb{Z}} = \min_{g, w} \sum_{n \in \hat{T}} \sum_{l \in L} \sum_{i \in I_l} (\mathcal{D}_{il} + 1)(1 - \hat{\gamma}_{il}^t)(\mathcal{g}_{lil}^{+n\mathcal{T}} + \mathcal{g}_{lil}^{-n\mathcal{T}} + \mathrm{w}_{lil}^{+n\mathcal{T}} + \mathrm{w}_{lil}^{-n\mathcal{T}}) \quad (41)$$

$$g_l^n - \hat{g}_{lil}^{n\mathcal{T}} = \mathcal{g}_{lil}^{+n\mathcal{T}} - \mathcal{g}_{lil}^{-n\mathcal{T}}, \quad \forall i \in I_l, l \in L, n \in \hat{T}, \quad (42)$$

$$w_l^n - \widehat{w}_{lil}^{n\mathcal{T}} = \mathrm{w}_{lil}^{+n\mathcal{T}} - \mathrm{w}_{lil}^{-n\mathcal{T}}, \quad \forall i \in I_l, l \in L, n \in \hat{T}, \quad (43)$$

$$\mathcal{g}_{lil}^{+n\mathcal{T}}, \mathcal{g}_{lil}^{-n\mathcal{T}}, \mathrm{w}_{lil}^{+n\mathcal{T}}, \mathrm{w}_{lil}^{-n\mathcal{T}} \geq 0, \quad \forall i \in I_l, l \in L, n \in \hat{T}. \quad (44)$$

The agreement subroutine for trajectories and signal timing plans is further detailed in Fig. 8.

| Agreement subroutine |
|---|
| **Step 1: Set** |
|     $\mathcal{T} \leftarrow 1$ |
|     $\hat{a}_{il}^{t\mathcal{T}} \leftarrow a_{il}^t \quad \forall i \in I_l, l \in L, t \in [t_0, t_0 + \widehat{N}/\Delta T]$ |
|     $\hat{g}_{lil}^{n\mathcal{T}} \leftarrow g_l^n, \widehat{w}_{lil}^{n\mathcal{T}} \leftarrow w_l^n \quad \forall l \in L, n \in [n_0, n_0 + \widehat{N}/\Delta \hat{T}]$ |
| **Step 2: Loop until agreements are reached** |
|     2.1. Optimization |
|         2.1.1. Solve problem $z$ for CAVs |
|         2.1.2. Solve problem $z'$ for CHVs |
|     2.2. Update |
|         2.2.1. Update CAV trajectories using equations (35) |
|         2.2.2. Update CHV trajectories using equations (36) |
|     2.3. Signal Agreement |
|         2.3.1. If signal parameters are fixed go to step 2.4. |
|         2.3.2. If votes on signals do not change solve problem (14)-(24), (38)-(39), and (41)-(44) and fix signal variables |
|         2.3.3. Set $\mathcal{T} \leftarrow \mathcal{T} + 1$ and go to step 2.1. |
|     2.4. Trajectory Agreement |
|         2.4.1. If $|\hat{x}_{il}^{t,\mathcal{T}+1} - \hat{x}_{il}^{t,\mathcal{T}}| \leq \varepsilon \quad \forall i \in I_l, l \in L, t \in T$ end |
|         2.4.2. Set $\mathcal{T} \leftarrow \mathcal{T} + 1$ and go to step 2.1. |
| **End Loop** |

**Fig. 8.** Agreement formation subroutine

## IV. NUMERICAL EXPERIMENTS

This section first introduces the case study intersection and then discusses the results including mobility, safety, and convergence performances. The proposed problem framework is run on a desktop computer with a Core i9 CPU and 64 GB of memory. Vissim microscopic traffic simulator [38] is used to implement the obtained solutions in order to measure mobility and safety performance measures. We overwrite Vissim's car following logic for CAVs based on the trajectories that our methodology finds. On the other hand, we let the Vissim car-following to update the trajectories of CHVs and use the car-following model (1) to only estimate the movements of CHVs.

### A. Case Study

The proposed problem methodology is applied to a four-legged isolated intersection with exclusive through and left-turn movements. We assume that vehicles are in their preferred lanes prior to entering the intersection neighborhood and thus, lane changing can be ignored. Please note that each phase is associated with one lane since the intersection has one lane for each movement. The communication range is assumed to be $650 \, ft$ which is within the maximum $1000 \, ft$ dedicated short-range communication range. The maximum speed for all vehicles is assumed to be $42.5 \, ft/s$. We update trajectories and signal indications every $\Delta T = 0.5 \, s$ and $\Delta \hat{T} = 2 \, s$, respectively. In addition, the maximum length of vehicle-groups must allow a vehicle group to clear the intersection with the minimum white phase time plus yellow transition time. Therefore, the maximum vehicle group length should be greater than $\bar{\zeta} \geq \bar{v}(\underline{W} + Y) + D = 351.8 \, ft$, which is rounded up to $360 \, ft$. TABLE II summarizes the value of the parameters used in the case studies.

TABLE II VALUE OF THE PARAMETERS USED IN THE CASE STUDY.

| Parameters | Value |
|---|---|
| Maximum and minimum speeds ($ft/s$) | 42.5 & 0 |
| Maximum acceleration ($ft/s^2$) | 13 |
| Minimum acceleration ($ft/s^2$) | -11.5 |
| Human driver reaction time ($s$) | 1 |
| Automated vehicle reaction time ($s$) | 0.1 |
| Detection range ($ft$) | 650 |
| Intersection location ($ft$) | 650 |
| Maximum length of vehicle-groups ($ft$) | 360 |
| Average length of vehicles ($ft$) | 13 |
| Safety distance between consecutive vehicles ($ft$) | 11.8 |
| Safety distance between CAVs and red signal light ($ft$) | 1 |
| Safety distance between vehicle-groups of conflicting lane-groups ($ft$) | 40 |
| Car-following parameter $\alpha_1$ ($s^{-1}$) | 0.95 |
| Car-following parameter $\alpha_2$ ($s^{-2}$) | 0.25 |
| Trajectory updating interval ($s$) | 0.5 |
| Signal status updating interval ($s$) | 2 |
| All red time ($s$) | 2 |
| Yellow time ($s$) | 4 |
| Minimum active (green + following white) time for through movements ($s$) | 12 |
| Minimum active (green + following white) time for left-turns ($s$) | 4 |
| Minimum white duration for through movements ($s$) | 6 |
| Minimum white duration for left-turning movements ($s$) | 4 |
| Maximum green time for through and left-turning movements ($s$) | 60 |
| Planning time horizon ($s$) | 20 |
| Study period ($s$) | 900 |

The problem is solved for a 900 $s$ study period and three different demand levels as shown in Fig. 9. Note that the demand for left-turn movements is 8% of the through movements. Eight different CAV market penetration rates are considered for each demand level, ranging from low (0% and 10), medium (30%, 50%, and 70%), to high (80%, 90%, and 100%) to capture the effects of CAVs on the intersection performance measures.

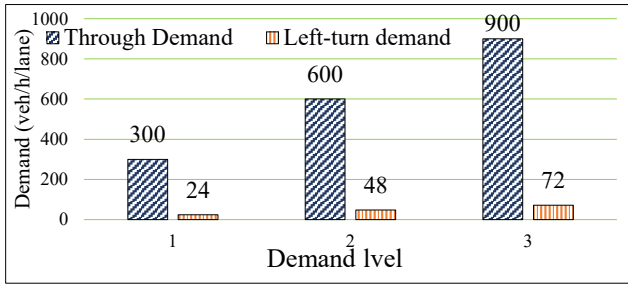

**Fig. 9.** Demand levels for case studies.

### B. Results

#### 3) Mobility

Fig. 10 shows the average delay for the fully-actuated signal timing plan found in Vissim [38] and the proposed agreement-based signal timing and trajectory optimization for different demand levels and CAV market penetration rates. The results indicate that the proposed methodology reduces the intersection average delay by 40.2% - 98.9% compared to the fully-actuated signal timing plan for various demand levels and CAV penetrations rates. Increasing the CAV penetration rate contributes to reductions in the average delay in all demand levels. This reduction can be caused by (i) more efficient operation of CAVs and (ii) increasing the possibility of activation of white phases. The significant reduction in intersection average delay from 80% to 90% CAV penetration rate in demand level 3 supports the hypothesis of the effectiveness of the white phase for a mixed autonomy environment. In addition, the difference between average delay for different demand levels decreases as the CAV market penetration rate increases. This can be due to efficient operation of vehicles during white phases which increases the intersection capacity significantly.

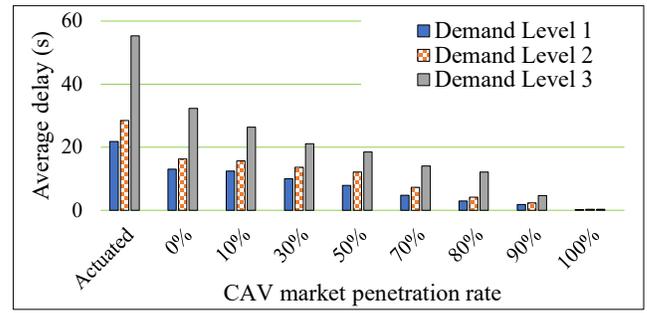

**Fig. 10.** Average delay for all demand levels with different CAV market penetration rates.

Fig. 11 shows the signal timing plans for all lane groups with different CAV market penetration rates in demand level 3. As shown in this figure, white phase activation increases with the CAV market penetration rate. We observe frequent white phase use for conflicting movements when the market penetration rate reaches 50%. The white phase becomes the dominant signal indication when the market share reaches and exceeds 70%; however, red and green signal indications appear even at a 90% CAV penetration rate.

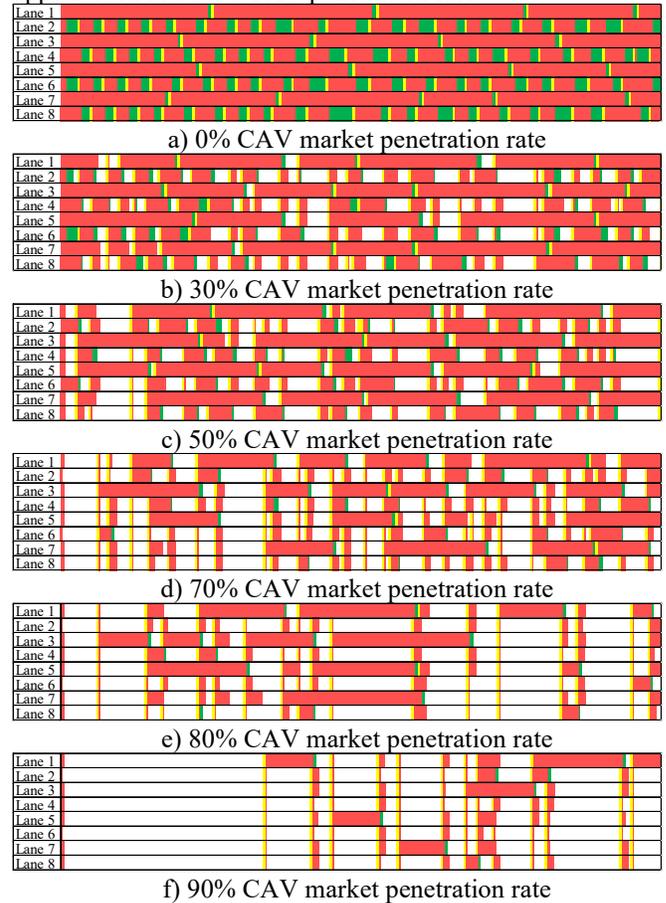

a) 0% CAV market penetration rate

b) 30% CAV market penetration rate

c) 50% CAV market penetration rate

d) 70% CAV market penetration rate

e) 80% CAV market penetration rate

f) 90% CAV market penetration rate

**Fig. 11.** Signal timing plans for demand level 3 with different CAV market penetration rates

A comparison between the activation rate of white phases with the central formulation proposed by Niroumand et al. [22] and the distributed formulation proposed in this study is presented in Fig. 12. One limitation of the previous work was that the white phase was not activated as frequently as it was expected. As such, in the present study, the formulation was enhanced. As it is shown in this figure, our proposed formulation leads to significantly higher rates of white phase activation than previous formulation for all CAV penetration rates due to the changes made to the formulation. The difference between the white phase activation rates affirms that our proposed formulation could successfully increase the feasible area of the problem.

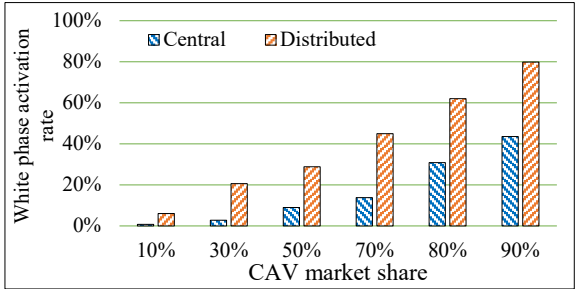

**Fig. 12.** White phase activation rate comparison between central and distributed formulations for demand level 3

TABLE III compares the average run time and total delay achieved by executing the central ([22], [30], [31]) and distributed methodologies for different CAV penetration rates under demand level 3. The average run time for the distributed methodology is always less than the updating time interval while the run time for the central methodology is much longer than the updating time interval by two orders of magnitude. Furthermore, the distributed methodology yielded shorter total delays than the central methodology when CAV penetration rate is between 50% and 90%. This is because of the higher white phase activation rate in the distributed formulation due to the enhancement we have made in the optimization model used in the central methodology. Note that with 0%, 10%, 30%, and 100% CAV penetration rates the signal timing plans cannot be significantly different between the central and distributed formulations since the possibility of activation of the white phases are low with 0%, 10%, and 30% CAV penetration rates while the white phase is the only phases used with 100% CAV penetration rates. Therefore, the total delay of the distributed methodology is more than the central one due to the decentralized methodology.

TABLE III  TOTAL DELAY AND AVERAGE RUN TIME COMPARISON BETWEEN CENTRAL AND DISTRIBUTED METHODOLOGIES UNDER DEMAND LEVEL 3

|   |   | Average run time (ms) | | Total delay (s) | |
|---|---|---|---|---|---|
|   |   | Distributed | Central | Distributed | Central |
| CAV penetration rate (%) | 0 | 167 | 10315 | 4905.7 | 4701.6 |
|  | 10 | 144 | 10649 | 3739 | 3604.4 |
|  | 30 | 165 | 10822 | 2644 | 2610.3 |
|  | 50 | 161 | 11034 | 2374 | 2436.1 |
|  | 70 | 104 | 13380 | 1602.6 | 1800 |
|  | 80 | 160 | 67684 | 1194 | 1356.3 |
|  | 90 | 128 | 115063 | 347.4 | 362.9 |
|  | 100 | 173 | 60930 | 85.2 | 82.31 |

Fig. 13 shows the trajectory of CHVs and CAVs on the eastbound lane of the intersection shown in Fig. 2 with different CAV penetration rates in demand level 3. The trajectory of CHVs and CAVs are shown with solid red and blue lines, respectively. Instead of stopping behind the red light, vehicles adjust their speed during white phases and pass the conflicting zone with a fewer number of stops. Therefore, traffic flow becomes smoother with an increase in the CAV market share.

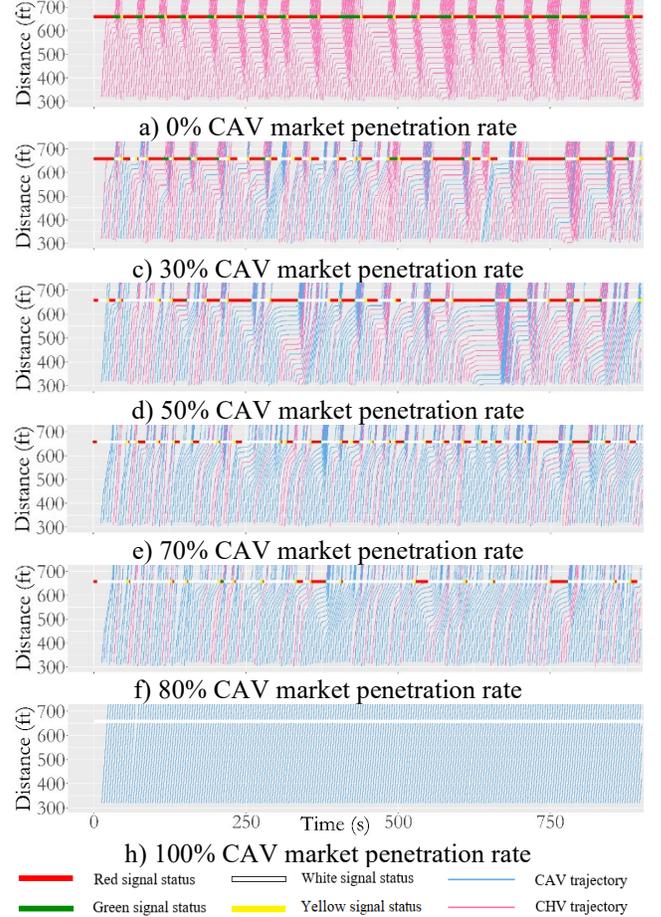

a) 0% CAV market penetration rate

c) 30% CAV market penetration rate

d) 50% CAV market penetration rate

e) 70% CAV market penetration rate

f) 80% CAV market penetration rate

h) 100% CAV market penetration rate

| Red signal status | White signal status | CAV trajectory |
| Green signal status | Yellow signal status | CHV trajectory |

**Fig. 13.** Trajectories of vehicles in lane 2 for demand level 3 with different CAV penetration rates

*4) Safety*

We have analyzed traffic safety considering rear-end and crossing conflicts using the surrogate safety assessment model (SSAM) software [39]. The time to collision (TTC) is used as the surrogate safety assessment measure. The results indicate that there is no near collision case for both considered conflicts with different CAV market penetration rates and demand levels assuming 1.5 $s$ TTC threshold. Fig. 14 shows the TTC for rear-end conflicts under the TTC threshold of 10 $s$. The median of TTC shows a decreasing trend by increasing the CAV penetration rate up to 80% due to co-existence of CHVs and CAVs with shorter reaction time. On the other hand, the median of the rear-end TTC decreases from 80% to 100% CAV penetration rates since most of the vehicles are operated during simultaneous white phases which leads to a smaller number of stops compared to green and red phases. On the other hand, interquartile dispersion increases as CAV penetration rate increases up to 80% and decreases from 80% to 100%. This can be due to mixed fleet of CAVs and CHVs and mixed use of green and white signal indications. However, the interquartile dispersion decreases as the fleet and signal indications become uniform in 90% and 100% CAV penetration rates. Note that, no crossing near-crash condition was reported by SSAM.

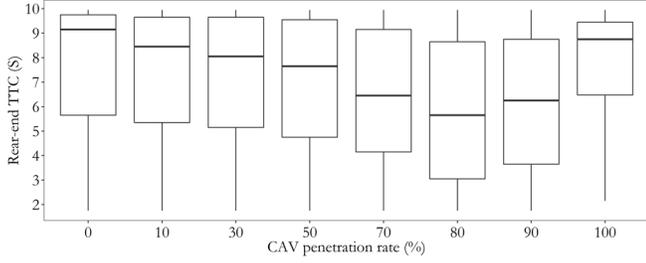

**Fig. 14.** Rear-end time to collision for demand level 3 with different CAV penetration rates

*5) Rider Comfort*

Fig. 15 shows the trajectory, speed, and acceleration profiles of a sample CAV and its following CHV that go through the intersection during the white phase. As it is shown, the CAV reduces its speed to prevent a potential collision with conflicting vehicle-groups. However, it does not need to come to a complete stop to avoid the collision. The following CHV (Fig. 15. (d), (e), and (f)) goes through the intersection even smoother than the leading CAV since the CAV is responsible for the collision avoidance.

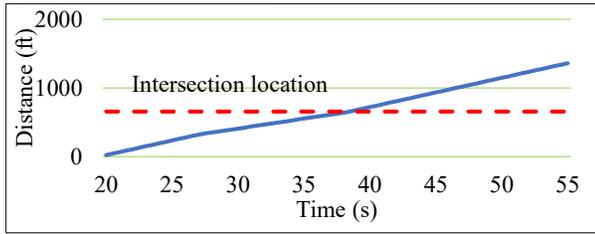

a) CAV trajectory

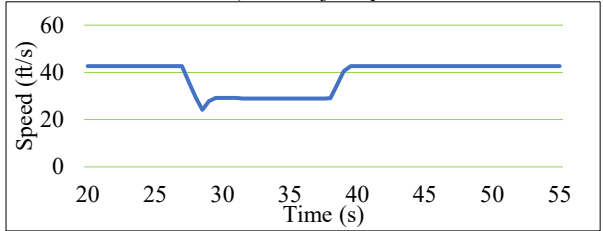

b) CAV speed profile

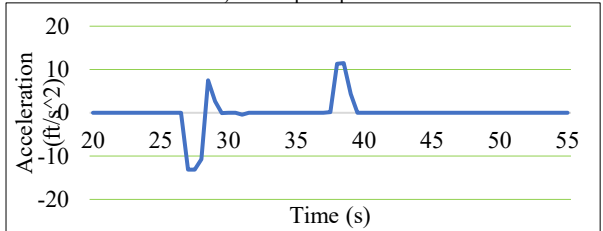

c) CAV acceleration profile

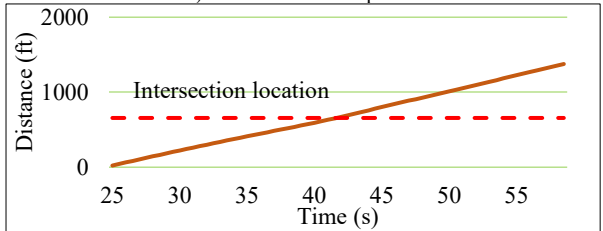

d) CHV trajectory

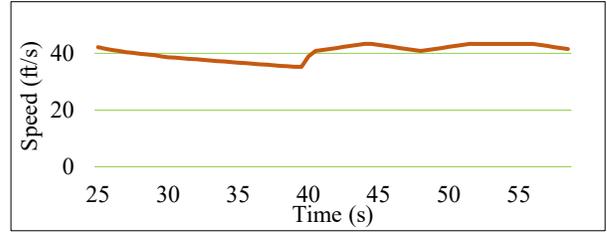

e) CHV speed profile

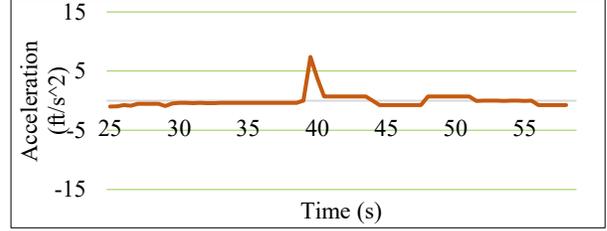

f) CHV acceleration profile

**Fig. 15.** Trajectory, speed, and acceleration profiles of a CAV and its following CHV

TABLE III summarizes the driving comfort parameters including average and standard deviation of speed, average acceleration and deceleration, average positive and negative jerk. Note that we did not include zero values in our analysis for the positive and negative jerk, acceleration, and deceleration. As shown in the table, the average speed increases for both types of vehicles as the CAV market penetrations rate increases. The average speeds for CAVs are higher than CHVs for all CAV penetration rates since CAVs have smaller reaction times and consider the future states of the signal and other CAVs in their trajectory optimization program. On the other hand, speed variances show a decreasing trend with increasing the CAV penetration rate since the speed difference of CAVs is minimized in their optimization program. Therefore, the speed difference of CAVs is always less than CHVs. In addition, average positive and negative acceleration and jerk are much lower than their maximum and minimum values which shows smooth changes in variations of speed and acceleration rate.

TABLE IV  DRIVING COMFORT PARAMETERS UNDER DEMAND LEVEL 3 AND DIFFERENT CAV PENETRATION RATES FOR CAVs AND CHVs

| Driving Comfort Measure | CAV Market Share (%) | | | | | |
|---|---|---|---|---|---|---|
| | 0 | 30 | 50 | 70 | 80 | 100 |
| CAVs | | | | | | |
| Average speed ($ft/s$) | - | 25.66 | 26.47 | 28.93 | 35.68 | 38.74 |
| Speed standard deviation ($ft/s$) | - | 18.19 | 17.87 | 16.69 | 12.01 | 5.87 |
| Average acceleration ($ft/s^2$) | - | 3.89 | 3.80 | 4.54 | 5.12 | 4.98 |
| Acceleration standard deviation ($ft/s^2$) | - | 4.23 | 4.12 | 4.49 | 4.66 | 4.29 |
| Average deceleration ($ft/s^2$) | - | -4.03 | -3.92 | -4.74 | -5.34 | -7.17 |
| Deceleration standard deviation ($ft/s^2$) | - | 5.07 | 4.93 | 5.10 | 5.29 | 5.85 |
| Average positive jerk ($ft/s^3$) | - | 3.89 | 3.98 | 5.15 | 7.55 | 8.31 |
| Positive jerk standard deviation ($ft/s^3$) | - | 5.99 | 5.91 | 6.91 | 8.62 | 8.14 |
| Average negative jerk ($ft/s^3$) | - | -3.83 | -4.39 | -5.20 | -7.66 | -10.90 |
| Negative jerk standard deviation ($ft/s^3$) | - | 6.04 | 6.18 | 7.13 | 8.82 | 9.20 |
| CHVs | | | | | | |
| Average speed ($ft/s$) | 21.72 | 24.39 | 24.97 | 25.71 | 31.39 | - |
| Speed standard deviation ($ft/s$) | 19.81 | 19.08 | 18.82 | 18.19 | 16.52 | - |
| Average acceleration ($ft/s^2$) | 3.39 | 2.99 | 3.07 | 2.84 | 1.94 | - |
| Acceleration standard deviation ($ft/s^2$) | 4.45 | 4.17 | 4.18 | 4.03 | 3.27 | - |
| Average deceleration ($ft/s^2$) | -3.39 | -2.45 | -2.15 | -1.95 | -1.29 | - |
| Deceleration standard deviation ($ft/s^2$) | 4.29 | 3.63 | 3.27 | 2.79 | 1.95 | - |
| Average positive jerk ($ft/s^3$) | 2.33 | 2.35 | 2.36 | 2.08 | 1.22 | - |
| Positive jerk standard deviation ($ft/s^3$) | 4.25 | 4.66 | 4.91 | 4.72 | 3.78 | - |
| Average negative jerk ($ft/s^3$) | -2.52 | -2.30 | -2.30 | -1.91 | -1.14 | - |
| Negative jerk standard deviation ($ft/s^3$) | 3.80 | 4.17 | 4.38 | 3.97 | 2.85 | - |

*6) White Phase Effect on Traffic Operations*

We study the effects of the white phase on total delay by creating two scenarios of (1) no white phase activation and (2) optimized white phase activation. The no white phase scenario represents joint CAV trajectory and signal optimization. In Fig. 16 we summarize the results for demand level 3 with 900 vehicle/h/lane on each through lane and 8% left turns. At 0% CAV market penetration rates, the total delays found for each scenario were identical as expected. In fact, the white phase is never activated in either scenario. At 10% CAV market share, we observed the white phase in 6.03% of the study period, which led to a 3.2% reduction in total delay due to reducing the number of phase transitions. With an increase in the CAV market share, the total delay reduction increases because, white phase allows much fewer phase transitions, which leads to less lost time and less delay.

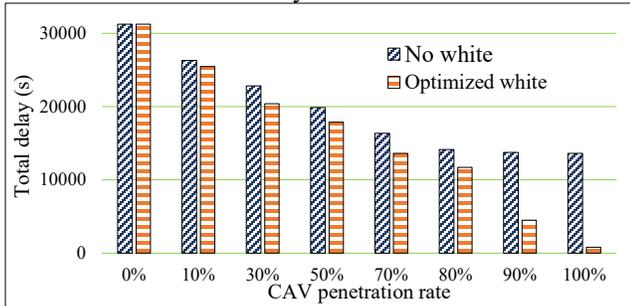

**Fig. 16.** No white phase and Optimized white phase total delay comparison for demand level 3

*7) Fuel Consumption*

Fig. 17 shows the average fuel consumption for different CAV penetration rates and white phase activation scenarios under demand level 3. As it is shown, average fuel consumption strictly decreases as the CAV penetration rate increases for both white phase activation scenarios. In addition, fuel consumption under optimized white phase activation scenario is always lower than no white phase activation scenario due to the presence of simultaneous white phases, which lead to fewer phase transitions and lost time. Average fuel consumption decreases with a relatively high rate from 0% to 30% CAV penetration rates due to higher performances of CAVs. The reduction rate decreases from 30% to 80% and further increases from 80% to 90% due to the dominance of the simultaneous white phases under optimized white phase activation scenario. On the other hand, the reduction in average fuel consumption under no white phase activation scenario is not significant in high CAV penetration rates since vehicles still need to decelerate to stop behind the red signal indications and accelerate to pass the conflicting area during green signals.

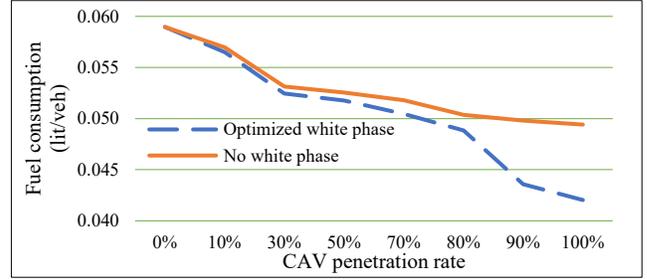

**Fig. 17.** No white phase and Optimized white phase average fuel consumption comparison for demand level 3

*8) Convergence*

The convergence of the proposed methodology is studied by comparing the trajectories of vehicles at different iterations. Fig. 18 shows the trajectory from 180th to 190th second. Note that this vehicle entered the detection area at the 170th second. The vehicle smoothly reduced its speed to enter the conflicting area at the proper time, however, it had to make a minor change in its trajectory prior to entering the conflicting area to avoid collisions. We also observe that the trajectory changes over iterations and the amount of change decreases over iterations until it becomes negligible. This is the point that the agreement is reached among CAVs on their trajectories. Fig. 19 shows the difference in trajectories between two consecutive iterations. The difference between iterations one and two is large; however, the difference decreases as the algorithm proceeds in iterations. Fig. 18 and Fig. 19 show that the trajectory of the vehicle is converged at iteration 7 with no considerable fluctuations. Note that, similar trends are observed for all the vehicles and we show this trajectory as a sample.

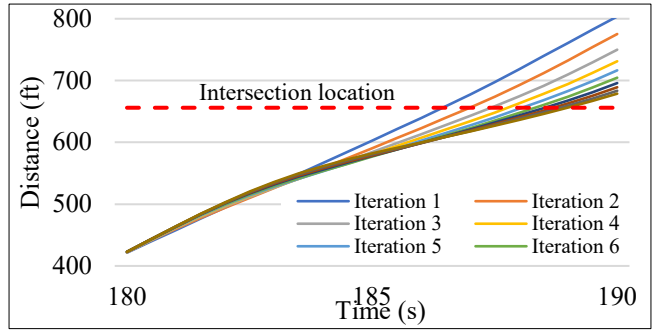

**Fig. 18.** Trajectory of a vehicle across different iterations

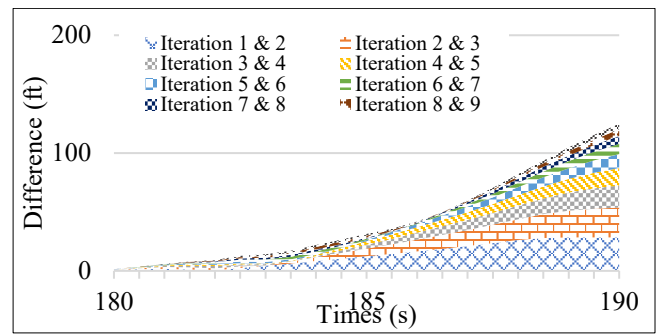

**Fig. 19.** The difference of trajectories at consecutive iterations

*9) Computation time*

Fig. 20 shows the computation time for demand level 3 with a 100% CAV penetration rate which is the most complicated scenario due to the presence of the highest number of vehicle-groups. As it is displayed, the algorithm converged with less than 0.4 seconds of computation time which represents the possibility of real-world implementation of the proposed methodology considering 0.5 seconds of updating time interval.

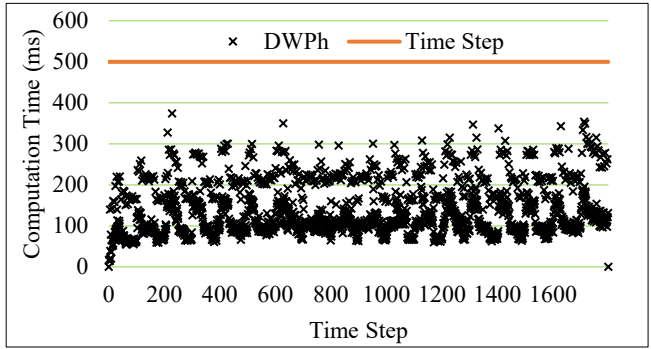

**Fig. 20.** Computation time for demand level 3 with 100% CAV penetration rate

## V. CONCLUSION

This paper introduces an agreement-based distributed methodology to control a mixed traffic stream of CAVs and CHVs at signalized intersections. The proposed solution technique aims to design CAV trajectories and signal timings more accurately while reducing the computational complexity of the decision-making process. We have incorporated a new white phase into traffic lights during which, CAVs act as mobile traffic controllers by forming and leading groups of CHVs and navigating through the intersection. The white phase will not be activated when the CAV market penetration rate is not high enough and instead, the green signal indication is utilized. We have formulated the joint signal timing and trajectory optimization problem as vehicle-level MINLPs. To reduce the complexity of the vehicle-level problems, the formulation is first linearized and then embedded into a receding horizon framework. The distributed methodology achieves agreement among all vehicles on vehicles' trajectories and signal timing parameters through an iterative process. Our case study results indicate that the proposed methodology can efficiently solve the problem, where the intersection total delay is reduced by 40.2% - 98.9% compared to a fully-actuated intersection control. In addition, the methodology produces solutions in real-time, which makes it suitable for real-world applications. Moreover, white phases are occasionally assigned to conflicting movements from a 30% CAV penetration rate and become the dominant signal status from a 70% rate.

## VI. REFERENCES

[1] R. Mohebifard and A. Hajbabaie, "Real-Time Adaptive Traffic Metering in a Connected Urban Street Network," in *Transportation Research Board 97th Annual Meeting*, 2018, no. 18–06061.
[2] R. Mohebifard and A. Hajbabaie, "Dynamic traffic metering in urban street networks: Formulation and solution algorithm," *Transp. Res. Part C Emerg. Technol.*, vol. 93, pp. 161–178, Aug. 2018.
[3] M. Tajalli, R. Niroumand, and A. Hajbabaie, "Distributed Cooperative Trajectory and Lane changing Optimization of Connected Automated Vehicles: Freeway Segments with Lane Drop."
[4] N. J. Goodall, B. L. Smith, and B. Park, "Traffic signal control with connected vehicles," *Transp. Res. Rec.*, vol. 2381, no. 1, pp. 65–72, 2013.
[5] W. Li and X. Ban, "Connected vehicles based traffic signal timing optimization," *IEEE Trans. Intell. Transp. Syst.*, 2018.
[6] K. Dresner and P. Stone, "A multiagent approach to autonomous intersection management," *J. Artif. Intell. Res.*, vol. 31, pp. 591–656, 2008.
[7] A. Mirheli, M. Tajalli, L. Hajibabai, and A. Hajbabaie, "A consensus-based distributed trajectory control in a signal-free intersection," *Transp. Res. part C Emerg. Technol.*, vol. 100, pp. 161–176, 2019.
[8] R. Mohebifard and A. Hajbabaie, "Trajectory Control in Roundabouts with a Mixed-fleet of Automated and Human-driven Vehicles," *Comput. Civ. Infrastruct. Eng.*, 2021.
[9] R. Mohebifard and A. Hajbabaie, "Connected Automated Vehicle Control in Single Lane Roundabouts," *Transp. Res. Part C Emerg. Technol.*, 2021.
[10] P. T. Li and X. Zhou, "Recasting and optimizing intersection automation as a connected-and-automated-vehicle (CAV) scheduling problem: A sequential branch-and-bound search approach in phase-time-traffic hypernetwork," *Transp. Res. Part B Methodol.*, vol. 105, pp. 479–506, 2017.
[11] Y. Guo, J. Ma, C. Xiong, X. Li, F. Zhou, and W. Hao, "Joint optimization of vehicle trajectories and intersection controllers with connected automated vehicles: Combined dynamic programming and shooting heuristic approach," *Transp. Res. part C Emerg. Technol.*, vol. 98, pp. 54–72, 2019.
[12] X. Liang, S. I. Guler, and V. V Gayah, "Joint optimization of signal phasing and timing and vehicle speed guidance in a connected and autonomous vehicle environment," *Transp. Res. Rec.*, vol. 2673, no. 4, pp. 70–83, 2019.
[13] M. Tajalli and A. Hajbabaie, "Distributed optimization and coordination algorithms for dynamic speed optimization of connected and autonomous vehicles in urban street networks," *Transp. Res. Part C Emerg. Technol.*, vol. 95, 2018.
[14] H. Rakha and R. K. Kamalanathsharma, "Eco-driving at signalized intersections using V2I communication," in *2011 14th international IEEE conference on intelligent transportation systems (ITSC)*, 2011, pp. 341–346.
[15] Z. Li, L. Elefteriadou, and S. Ranka, "Signal control optimization for automated vehicles at isolated signalized intersections," *Transp. Res. Part C Emerg. Technol.*, vol. 49, pp. 1–18, 2014.
[16] C. Yu, Y. Feng, H. X. Liu, W. Ma, and X. Yang, "Integrated optimization of traffic signals and vehicle trajectories at isolated urban intersections," *Transp. Res. Part B Methodol.*, vol. 112, pp. 89–112, 2018.
[17] K. Yang, S. I. Guler, and M. Menendez, "Isolated intersection control for various levels of vehicle technology: Conventional, connected, and automated vehicles," *Transp. Res. Part C Emerg. Technol.*, vol. 72, pp. 109–129, 2016.
[18] H. Jung, S. Choi, B. B. Park, H. Lee, and S. H. Son, "Bi-level optimization for eco-traffic signal system," in *2016 International Conference on Connected Vehicles and Expo (ICCVE)*, 2016, pp. 29–35.
[19] B. Xu et al., "Cooperative method of traffic signal optimization and speed control of connected vehicles at isolated intersections," *IEEE Trans. Intell. Transp. Syst.*, vol. 20, no. 4, pp. 1390–1403, 2018.
[20] W. Zhao, D. Ngoduy, S. Shepherd, R. Liu, and M. Papageorgiou, "A platoon based cooperative eco-driving model for mixed automated and human-driven vehicles at a signalised intersection," *Transp. Res. Part C Emerg. Technol.*, vol. 95, pp. 802–821, 2018.
[21] C. Chen, J. Wang, Q. Xu, J. Wang, and K. Li, "Mixed platoon control of automated and human-driven vehicles at a signalized intersection: dynamical analysis and optimal control," *Transp. Res. Part C Emerg. Technol.*, vol. 127, p. 103138, 2021.
[22] R. Niroumand, M. Tajalli, L. Hajibabai, and A. Hajbabaie, "Joint Optimization of Vehicle-group Trajectory and Signal Timing: Introducing the White Phase for Mixed-autonomy Traffic Stream,"

[23] X. Qian, J. Gregoire, F. Moutarde, and A. De La Fortelle, "Priority-based coordination of autonomous and legacy vehicles at intersection," in *17th international IEEE conference on intelligent transportation systems (ITSC)*, 2014, pp. 1166–1171.

[24] M. Li, X. Wu, X. He, G. Yu, and Y. Wang, "An eco-driving system for electric vehicles with signal control under V2X environment," *Transp. Res. Part C Emerg. Technol.*, vol. 93, pp. 335–350, 2018.

[25] M. Pourmehrab, L. Elefteriadou, S. Ranka, and M. Martin-Gasulla, "Optimizing signalized intersections performance under conventional and automated vehicles traffic," *IEEE Trans. Intell. Transp. Syst.*, 2019.

[26] H. Moradi, S. Sasaninejad, S. Wittevrongel, and J. Walraevens, "Proposal of an integrated platoon-based Round-Robin algorithm with priorities for intersections with mixed traffic flows," *IET Intell. Transp. Syst.*, vol. 15, no. 9, pp. 1106–1118, 2021.

[27] Y. Du, W. ShangGuan, and L. Chai, "A coupled vehicle-signal control method at signalized intersections in mixed traffic environment," *IEEE Trans. Veh. Technol.*, vol. 70, no. 3, pp. 2089–2100, 2021.

[28] M. Tajalli and A. Hajbabaie, "Traffic Signal Timing and Trajectory Optimization in a Mixed Autonomy Traffic Stream," *IEEE Trans. Intell. Transp. Syst.*, 2021.

[29] D. Rey and M. W. Levin, "Blue phase: Optimal network traffic control for legacy and autonomous vehicles," *Transp. Res. part B Methodol.*, vol. 130, pp. 105–129, 2019.

[30] R. Niroumand, M. Tajalli, L. Hajibabai, and A. Hajbabaie, "The Effects of the 'White Phase' on Intersection Performance with Mixed-Autonomy Traffic Stream," in *2020 IEEE 23rd International Conference on Intelligent Transportation Systems (ITSC)*, 2020, pp. 1–6.

[31] R. Niroumand, L. Hajibabai, and A. Hajbabaie, "The Effects of Connectivity on Intersection Operations with 'White Phase,'" in *2021 IEEE International Intelligent Transportation Systems Conference (ITSC)*, 2021, pp. 3839–3844.

[32] R. Niroumand, L. Hajibabai, A. Hajbabaie, and M. Tajalli, "Effects of Autonomous Driving Behavior on Intersection Performance and Safety in the Presence of White Phase for Mixed-Autonomy Traffic Stream," *Transp. Res. Rec.*, p. 03611981221082580, 2022.

[33] W. Helly, "Simulation of bottlenecks in single-lane traffic flow," 1959.

[34] M. Mehrabipour and A. Hajbabaie, "A cell-based distributed-coordinated approach for network-level signal timing optimization," *Comput. Civ. Infrastruct. Eng.*, vol. 32, no. 7, pp. 599–616, 2017.

[35] M. Mehrabipour and A. Hajbabaie, "A cell based distributed-coordinated approach for network level signal timing optimization," *Comput. Aided Civ. Infrastruct. Eng. an Int. J.*, vol. 32, no. 7, pp. 599--616, 2017.

[36] S. M. A. B. Al Islam, A. Hajbabaie, and H. M. A. Aziz, "A real-time network-level traffic signal control methodology with partial connected vehicle information," *Transp. Res. Part C Emerg. Technol.*, vol. 121, p. 102830, Dec. 2020.

[37] W. B. Powell, *Approximate Dynamic Programming: Solving the curses of dimensionality*, vol. 703. John Wiley & Sons, 2007.

[38] PTV Group, "PTV Vissim 7 User Manual," 2015.

[39] D. Gettman and L. Head, "Surrogate safety measures from traffic simulation models," *Transp. Res. Rec.*, vol. 1840, no. 1, pp. 104–115, 2003.

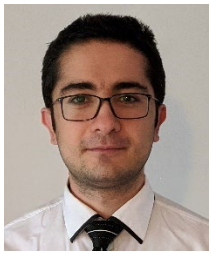
Ramin Niroumand received the B.Sc. degree in civil engineering from the University of Tabriz, Tabriz, Iran in 2014, and the M.Sc. degree in civil engineering from the Sharif University of Technology, Tehran, Iran, in 2017. He is currently pursuing the Ph.D. degree with the Department of Civil, Construction, and Environmental Engineering, North Carolina State University. His research interests include network optimization, traffic assignment, traffic signal optimization, and CAV trajectory optimization.

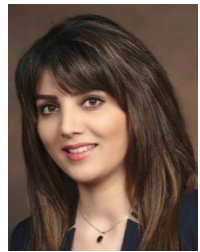
Leila Hajibabai is an assistant professor in the Department of Industrial and Systems Engineering at North Carolina State University (NCSU). She has received her Ph.D. from the University of Illinois at Urbana – Champaign (UIUC) in 2014 and has obtained two M.Sc. degrees, one from University of Tehran and the other from UIUC in Industrial Engineering and the other in Civil Engineering. Her research is in the transportation and logistics systems area. Dr. Hajibabai is a member of INFORMS, IEEE, and the Transportation Research Board (TRB)'s committees on Network Modeling, and Maintenance and Fleet Equipment. She is a chair of TRB's Young Members subcommittee of Highway Maintenance Section.

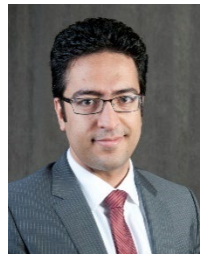
Ali Hajbabaie received the B.S. and M.S. degrees from Sharif University of Technology, Tehran, Iran in 2003 and 2006, respectively. He also received the M.S. degree in industrial engineering and the Ph.D. degree in civil engineering from the University of Illinois at Urbana-Champaign, in 2011 and 2012, respectively. He is an Associate Professor in the Civil, Construction, and Environmental Engineering Department at North Carolina State University. He is the recipient of the best paper award from the work zone traffic control committee of the transportation research board. He is the secretary of the work zone traffic control and the chair of the asset management subcommittee of the traffic signal systems committee of the transportation research board. His research interests include cooperative congestion management, distributed optimization, traffic operations, and traffic flow theory.